\documentclass[11pt,a4paper]{article}
\usepackage{bbm}

\usepackage[leqno]{amsmath}
\usepackage{amsfonts}
\usepackage{graphicx}
\usepackage{amsmath}
\usepackage{amssymb}
\usepackage{latexsym}
\usepackage{amsmath, amsfonts,amssymb, amsthm, euscript,makeidx,color,mathrsfs}
\usepackage{enumerate}
\usepackage{tabu}

\usepackage[colorlinks,linkcolor=blue,anchorcolor=green,citecolor=red]{hyperref}

\def\sqr#1#2{{\vcenter{\vbox{\hrule height.#2pt
              \hbox{\vrule width.#2pt height#1pt \kern#1pt \vrule width.#2pt}
              \hrule height.#2pt}}}}
\def\5n{\negthinspace \negthinspace \negthinspace \negthinspace \negthinspace }
\def\4n{\negthinspace \negthinspace \negthinspace \negthinspace }
\def\3n{\negthinspace \negthinspace \negthinspace }
\def\2n{\negthinspace \negthinspace }
\def\1n{\negthinspace }

\def\={\buildrel \triangle \over =}

\def\bal{\begin{aligned}}
\def\eal{\end{aligned}}
\def\no{\noindent}

\def\bs{\bigskip}

\def\liminf{\mathop{\underline{\rm lim}}}

\def\dim{\hbox{\rm dim$\,$}}

\def\sgn{\hbox{\rm sgn$\,$}}

\def\({\Big (}
\def\){\Big )}
\def\[{\Big[}
\def\]{\Big]}

\def\bde{\begin{definition}\label}
\def\ede{\end{definition}}
\def\be{\begin{equation}}
\def\bel{\begin{equation}\label}
\def\ee{\end{equation}}
\def\beq{\begin{equation*}\begin{aligned}}
\def\eeq{\end{aligned}\end{equation*}}
\def\bt{\begin{theorem}\label}
\def\et{\end{theorem}}
\def\bc{\begin{corollary}\label}
\def\ec{\end{corollary}}
\def\bl{\begin{lemma}\label}
\def\el{\end{lemma}}
\def\bp{\begin{proposition}\label}
\def\ep{\end{proposition}}
\def\bas{\begin{assumption}\label}
\def\eas{\end{assumption}}
\def\br{\begin{remark}\label}
\def\er{\end{remark}}
\def\bex{\begin{example}\label}
\def\ex{\end{example}}
\def\ba{\begin{array}}
\def\ea{\end{array}}
\def\ed{\end{document}}

\def\square#1{\vbox{\hrule\hbox{\vrule height#1%
     \kern#1\vrule}\hrule}}
\def\rectangle#1#2{\vbox{\hrule\hbox{\vrule height#1%
     \kern#2\vrule}\hrule}}

\font\tenbb=msbm10 \font\sevenbb=msbm7 \font\fivebb=msbm5

\newfam\bbfam
\scriptscriptfont\bbfam=\fivebb \textfont\bbfam=\tenbb
\scriptfont\bbfam=\sevenbb

\newlength{\defbaselineskip}
\newcommand{\setlinespacing}[1]%
           {\setlength{\baselineskip}{#1\defbaselineskip}}

\theoremstyle{plain}
\newtheorem{thm}{Theorem}[section]
\newtheorem{cor}[thm]{Corollary}
\newtheorem{lem}[thm]{Lemma}
\newtheorem{prop}[thm]{Proposition}
\newtheorem{exam}[thm]{Example}
\newtheorem{rem}[thm]{Remark}
\newtheorem{defi}[thm]{Definition}

\def\poly{\mathbb{C}[z_1,\ldots,z_d]}
\def\Hol{\mathrm{Hol~}}
\def\Id{\mathrm{Id}}

\def\card{\mathrm{card~}}

\def\sgn{\mathrm{sgn}}

\def\tensorspace{\bigotimes\limits_{i=1}^mH_i^m}

\makeatletter
   
   \@addtoreset{equation}{section}
\makeatother

\begin{document}
\title{\bf Essentially normal   quotient weighted  Bergman modules over the bidisk  and distinguished varieties }

\author{Kunyu Guo, Penghui Wang, Chong Zhao}

\date{Feb 19, 2023}
\maketitle

\begin{abstract}

We introduce a Grassmannian structure for a class of  quotient Hilbert modules and attack the polydisc version of Arveson-Douglas conjecture  associated to distinguished varieties. More interestingly,  we obtain an  operator-theoretic characterization of distinguished varieties  in the bidisk in terms of essential normality of the   quotient modules. As an application, we study the K-homology of the boundary of distinguished variety.

\end{abstract}

\bs


\bs

\no{\bf Key Words:} distinguished variety, essential normality, weighted Bergman module, Grassmannian quotient module, K-homology, algebraic variety.

\section{Introduction}

Algebraic zero variety is one of the main objects in algebraic geometry.  For an ideal $I$ of the polynomial ring $\mathbb C[z_1,\cdots,z_d]$, its algebraic zero variety is defined by
$$Z(I)=\{(z_1,\ldots,z_d)\in\mathbb{C}^d:p(z_1,\ldots,z_d)=0,~\forall p\in I\}.$$
A subset $V$ of the unit polydisc $\mathbb D^d$ is called an algebraic subvariety of $\mathbb D^d$, if there is an ideal $I$ of the polynomial ring, such that $V=Z(I)\cap\mathbb D^d$. A special subvariety, called the distinguished variety, was introduced by Agler and M\raisebox{.5ex}{c}Carthy \cite{AM1}, which plays an important role in Ando's inequality and Pick interpolation problem, etc. There are many efforts to explore the structure and applications of distinguished variety, please refer to  \cite{AM1,AM2,BKS,DGH,DS,DSS,Kne1,Kne2,PS,Sch,Ve} and references therein.
   \begin{defi}(Agler and M\raisebox{.5ex}{c}Carthy \cite{AM1})  Let $\mathbb{T}$ be the unit circle.
  A subvariety $V$ in $\mathbb D^2$ is called a distinguished variety, if $
  \partial V\subset\mathbb T^2.$
  \end{defi}

  Agler and M\raisebox{.5ex}{c}Carthy characterized distinguished varieties in terms of matrix-valued inner functions, where a matrix-valued function on $\mathbb D$ is called inner if it is analytic on $\mathbb D$ and unitary almost everywhere on $\mathbb T$.

\begin{thm}(Agler and M\raisebox{.5ex}{c}Carthy \cite{AM1})
    Let $V$ be a distinguished variety, then there is a rational matrix-valued inner function $\Psi$ such that
    $$V=Z\big(\det(wI_m-\Psi(z))\big)\cap\mathbb{D}^2.$$
\end{thm}
Recently, it is interesting to find that distinguished variety is closely related to essential normality of quotient Hilbert modules.
Recall that given a tuple $\underline{T}=(T_1,\ldots,T_d)$ of commuting operators on a Hilbert space $H$, one can naturally make $H$ into a Hilbert module \cite{DP} over the polynomial ring $\mathbb C[z_1,\ldots,z_d]$, with the module action
$$p\cdot x=p(T_1,\ldots, T_d)x,~p\in \mathbb C[z_1,\ldots,z_d],~x\in H.$$
For example, most of  the usual analytic Hilbert spaces over bounded domains \cite{CG}  are naturally  Hilbert modules  on $\poly$.
Motivated by the BDF-theory \cite{BDF}, one of the fundamental problems in the Hilbert module theory is to study the essential commutativity of the $C^*$-algebra $C^*(\underline{T})$ generated by $\{\Id_H,T_1,\ldots,T_d\}$ and the ideal $\mathcal{K}$ of compact operators, where $\Id_H$ is the identity on $H$, and the essential commutativity of two operators $A$ and $B$ means that the commutator $[A,B]=AB-BA$ is compact. When this happens, $H$ is said to be essentially normal.

A closed subspace $\cal M$ of $H$ is called a submodule if it is invariant under the module actions. Let  $\mathcal{N}=H\ominus\mathcal{M}$, and for $p\in\poly$ write
$$S_p=P_\mathcal{N}~p(T_1,\ldots, T_d)\mid_\mathcal{N},$$
the compression of $p(T_1,\ldots, T_d)$ on $\mathcal{N}$. Then $\mathcal{N}$ is naturally equipped with a $\mathbb C[z_1,\ldots,z_d]$-module structure via the tuple $(S_{z_1},\ldots, S_{z_d})$, and is called a quotient module of $H$.

In \cite{Ar1}, Arveson conjectured that graded submodules of the Drury-Arveson module $H_d^2$ over the unit ball in $\mathbb C^d$ are essentially normal. 
Much work has been done along this line, such as \cite{Ar2,Dou1,Dou2,DGW,DTY,DW,DWY,Guo,GWk1,GWk2,GW,KS,Sha,WX1,WX2,WX3}.

However, the  Arveson's conjecture is well-known to  fail with nonzero Hardy submodules over the polydisc and we shall turn to quotient  modules. As usual, the Hardy space over $\mathbb{D}^d$ is defined as
$$H^2(\mathbb{D}^d)=\left\{f\in\Hol(\mathbb{D}^d):||f||^2=\sup_{r\in[0,1)}\int_{\mathbb{T}^d}|f(rz)|^2\mathrm{d}m(z)<\infty\right\},$$
where $\mathrm{d}m$ is the normalized Haar measure on the torus $\mathbb{T}^d$.
 In \cite{GWp1,GWp2,WZ1,WZ3},   the essential normality of (quasi-)homogenous quotient modules was completely characterized. As  part of the results, it was proved that, for a (quasi-)homogenous ideal $I$ of $\mathbb C[z_1,\cdots,z_d]$, if $Z(I)\cap \partial \mathbb D^d\subset\mathbb T^d$, then the associated quotient module $[I]^\perp$ in $H^2(\mathbb D^d)$ is essentially normal.

For a first glance at the inhomogeneous case, let us recall the results in \cite{Cla,Wa}. Let $\eta_1$ and $\eta_2$ be two nonconstant inner functions on $\mathbb D$, then the  quotient module  $[\eta_1(z)-\eta_2(w)]^\perp$ of $H^2(\mathbb D^2)$ is essentially normal if and only if both $\eta_1$ and $\eta_2$ are finite Blaschke products, which is equivalent to the condition that $Z(\eta_1-\eta_2)\cap\mathbb{D}^2$ is a distinguished variety.

In the present paper, we will consider the essential normality of quotient modules in the weighted Bergman modules. For real number $\alpha>-1$, the weighted Bergman space is defined by
$$\mathcal{A}_\alpha(\mathbb{D}^d)=\left\{f\in\Hol(\mathbb{D}^d):||f||_\alpha^2=c_\alpha\int_{\mathbb{D}^d}|f(z)|^2\prod_{i=1}^d(1-|z_i|^2)^\alpha\mathrm{dV}_d(z)<\infty\right\},$$
where $\mathrm{dV}_d$ is the usual Lebesgue measure on $\mathbb{D}^d$, and $c_\alpha$ is the constant that normalizes the measure $\prod\limits_{i=1}^d(1-|z_i|^2)^\alpha\mathrm{dV}_d$. It is known that the weighted Bergman space is a reproducing kernel Hilbert space with reproducing kernel
\begin{equation}\label{eq:REP}
    K_{\lambda}(z)=\prod\limits_{i=1}^d{ (1-\bar{\lambda}_i z_i)^{-\alpha-2}},\lambda,z\in\mathbb{D}^d.
\end{equation}
Notice that when $\alpha=-1$, (\ref{eq:REP}) is just the reproducing kernel of the Hardy space $H^2(\mathbb{D}^d)$. In this sense, we also use $\mathcal{A}_{-1}(\mathbb{D}^d)$ to refer to $H^2(\mathbb D^d)$.

In what follows, for an ideal $I\subset\mathbb{C}[z,w]$, if $Z(I)\cap\mathbb{D}^2$ is a distinguished variety, we call $I$ a distinguished ideal, and both $[I]_\alpha\subset\mathcal{A}_\alpha(\mathbb{D}^2)$ the submodule generated by $I$ and $[I]_\alpha^\perp$ the associated quotient module are said to be distinguished. To simplify the notation, for $f\in\mathcal{A}_\alpha(\mathbb{D}^d)$, $[f]_\alpha$ denotes the principal submodule of $\mathcal{A}_\alpha(\mathbb{D}^d)$ generated by $f$. In particular, we shall write $[f]=[f]_{-1}$ in short, which is the submodule of Hardy module $H^2(\mathbb D^d)$.


In \cite{WZ2,WZ3}, it was proved that under the assumption of quasi-homogeneity, all the distinguished quotient modules of $\mathcal{A}_\alpha(\mathbb{D}^d)$ are essentlally normal. In fact, they constitute the major parts of essentially normal quotient modules. It is natural to ask whether all the distinguished quotient modules are essentially normal, without the quasi-homogeneity.

To state our main results, we introduce the notion of variety rank and Macaev interpolation ideal $\mathcal{L}^{(1,\infty)}$.
For a given distinguished variety $V$, set
$$m_V=\max_{\lambda\in\mathbb{D}}\card\{w\in\mathbb{D}:(\lambda,w)\in V\},~n_V=\max_{\lambda\in\mathbb{D}}\card\{z\in\mathbb{D}:(z,\lambda)\in V\}.$$
In the language of Agler and M\raisebox{.5ex}{c}Carthy \cite{AM1}, $(m_V,n_V)$ is called the rank of the variety $V$. $V$'s rank plays a critical role in our research on weighted Bergman modules.

For a compact operator $A\in B(H)$, we list the eigenvalues of $|A|=\sqrt{A^*A}$ in the decreasing order as $\{\lambda_0,\lambda_1,\ldots\}$, multiplicity being counted. The Macaev ideal $\mathfrak{L}^{(1,\infty)}$ consists of the compact operators for which the sequence $\Big\{{1\over \ln(n+1)}\sum\limits_{k=0}^n\lambda_k\Big\}$ is bounded. For a Hilbert $\poly$-module $H$, if $[p(\underline{T})^*,q(\underline{T})]\in\mathfrak{L}^{(1,\infty)}$ for all polynomials $p,q\in\mathbb C[z_1,\cdots,z_d]$, then $H$ is called $(1,\infty)$-essentially normal. Readers are referred to \cite{Conn} for more details of Macaev ideals.

The following theorem is the main result in this paper, which allows us to study distinguished variety from the viewpoint of Hilbert modules.
\begin{thm}\label{main-theorem}
Let $V$ be an algebraic subvariety of $\mathbb D^2$, then $V$ is distinguished if and only if there exists an integer $\alpha\geq\max\{m_{V},n_{V}\}-2,$ such that for every ideal $I$ with $V=Z(I)\cap \mathbb D^2,$ $[I]_{\alpha}^\perp$ is essentially normal. Furthermore, let $I\subset\mathbb{C}[z,w]$ be a distinguished ideal, then for every integer $$\alpha\geq\max\{m_{Z(I)},n_{Z(I)}\}-2,$$ the quotient module $[I]_\alpha^\perp$ of $\mathcal{A}_\alpha(\mathbb D^2)$ is $(1,\infty)$-essentially normal.
\end{thm}
For simplicity, we shall write $S_p$ for $S_p^\mathcal{N}$ in short, if there is no confusion. By the BDF-theory \cite{BDF}, when $[I]_\alpha^\perp$ is essentially normal, we have the following short exact sequence
\begin{eqnarray}\label{exact-sequence}
0\to {\cal K}\to C^*([I]_\alpha^\perp)\to C(\sigma_e(S_z,S_w))\to 0,
\end{eqnarray}
where $C^*([I]_\alpha^\perp)$ is the $C^*$-algebra generated by $S_z,~S_w$ and $\Id_{[I]_\alpha^\perp}$, and $\sigma_e(S_z,S_w)$ is the Taylor's joint essential spectrum \cite{Cur,Tay1} of the pair $(S_z,S_w)$. It will be shown that for distinguished ideals $I$,
$$\sigma_e(S_z,S_w)=\partial (Z(I)\cap \mathbb D^2).$$
In this paper, if $Z(I)\cap \mathbb D^2$ is a distinguished subvariety, it is proved that the short exact sequence (\ref{exact-sequence}) is not splitting, which gives a nontrivial element $e_I$ in the K-homology group $K\big(\partial (Z(I)\cap \mathbb D^2)\big)$.

The present paper is organized as follows. In Section 2, we prove the injectivity of both $S_z$ and $S_w$, which reduces the essential normality of $[I]_\alpha^\perp$ to the essential co-isometricity of $S_z$ and $S_w$. In Section 3, the construction of Grassmannian quotient module is introduced, which is a new tool to represent the determinant of module-action-valued matrices. In Section 4, we prove the main result Theorem 1.3. In Section 5, the essential Taylor spectrum of the quotient modules and the related $K$-homology are studied. In Section 6, we discuss the decomposition of essentially normal quotient modules, in terms of the ones with lower weight index. In Section 7, we show some non-distinguished examples of essentially normal quotient modules.

\section{Essential isometry and essential co-isometry}\label{subsection:main}

~~~~This section is devoted to reducing the proof of the essential normality of distinguished quotient module $[I]_\alpha^\perp$ of $\mathcal{A}_\alpha(\mathbb{D}^2)$ to the compactness of $\Id_{[I]_\alpha^\perp}-S_zS_z^*$ and $\Id_{[I]_\alpha^\perp}-S_wS_w^*$. To begin with, we make the following reduction to the ideals in consideration.
\begin{lem}\label{lem:contracted}
    Let $I\subset\mathbb{C}[z,w]$ be an ideal, and $\alpha\geq-1$, then there is a unique ideal $I_0=[I]_\alpha\cap\mathbb{C}[z,w]$ such that $[I]_\alpha=[I_0]_\alpha$ and each algebraic component of $Z(I_0)$ intersects $\mathbb{D}^2$. If $I$ is distinguished, then $I_0$ is also distinguished.
\end{lem}
\begin{proof}
    If each algebraic component of $Z(I)$ intersects $\mathbb{D}^2$, the existence of $I_0$ is obvious. So in the remaining we assume at least one algebraic component of $Z(I)$ does not intersect $\mathbb{D}^2$.

    Let $I=\bigcap_{i=1}^k I_i$ be the primary decomposition, and let $I_0$ be the intersection of those $I_i$ whose variety intersects $\mathbb{D}^2$, then each algebraic component of $Z(I_0)$ intersects $\mathbb{D}^2$. If $I$ is assumed distinguished, then so is $I_0$. From $I\subset I_0$ it follows $[I]_\alpha\subset[I_0]_\alpha$, and we only need to prove the inverse inclusion. If some $Z(I_i)$ does not intersect $\mathbb{D}^2$, then we can choose $p_i\in\sqrt{I_i}$ such that $Z(p_i)$ is not contained in $Z(I_0)$, and $Z(p_i)$ does not intersect $\mathbb{D}^2$. Suppose $p_i^{d_i}\in I_i$, and set $p=\prod_{Z(I_i)\cap\mathbb{D}^2=\emptyset}p_i^{d_i}$, then $pI_0\subset I$ and $Z(p)$ does not intersect $\mathbb{D}^2$. By \cite{Ge1,Ge2} or \cite[Theorem 2.2.13]{CG}, $p\mathbb{C}[z,w]$ is dense in both $H^2(\mathbb{D}^2)$ and $L_a^2(\mathbb{D}^2)$. By a careful checking, the proofs for weighted Bergman spaces remain valid. Therefore
    $$[I_0]_\alpha=\overline{I_0\mathcal{A}_\alpha(\mathbb{D}^2)}=\overline{I_0p\mathbb{C}[z,w]}=[pI_0]_\alpha\subset[I]_\alpha\subset[I_0]_\alpha,$$
    which gives $[I_0]_\alpha=[I]_\alpha$.

    The uniqueness of $I_0$ comes from \cite{DP}. Since each component of $Z(I_0)$ intersects $\mathbb{D}^2$, \cite{DP} ensures $$I_0=[I_0]_\alpha\cap\mathbb{C}[z,w]=[I]_\alpha\cap\mathbb{C}[z,w].$$
\end{proof}
\begin{rem}
By \cite{DP}, an ideal $I\subset\mathbb{C}[z,w]$ is called contracted if $[I]_\alpha\cap\mathbb{C}[z,w]=I$. By \cite{CG,Ge1,Ge2}, $I$ is contracted if and only if each component of its variety intersects $\mathbb{D}^2$. According to Lemma \ref{lem:contracted}, every distinguished submodule of $\mathcal{A}_\alpha(\mathbb{D}^2)$ is generalized by a contracted distinguished ideal. For this reason, if there is no special statement, throughout this paper all the distinguished ideals mentioned are assumed to be contracted.
\end{rem}
The proof of the main theorem is based on observation that the essential normality of the distinguished quotient module $[I]_\alpha^\perp$ can be reduced to the compactness of both $\Id_{[I]_\alpha^\perp}-S_zS_z^*$ and $\Id_{[I]_\alpha^\perp}-S_wS_w^*$. To see this, we need the Fredholmness of $S_z$ and $S_w$.
\begin{lem}\label{lem:injectivity1}
    Let $I\subset\mathbb{C}[z,w]$ be an ideal, $a\in\mathbb{C}\backslash\mathbb{T}$ such that $Z(I)\cap(\{a\}\times\mathbb{T})=\emptyset$, then $\ker S_{z-a}$ is finite dimensional. Symmetrically, if $Z(I)\cap(\mathbb{T}\times\{a\})=\emptyset$, then $\ker S_{w-a}$ is finite dimensional. 
\end{lem}
\begin{proof}
    By symmetry, we only prove the conclusions for $S_{z-a}$. The conclusions are obvious for $|a|>1$, so we assume $|a|<1$ in the remaining of the proof.

    First we assume that $I$ is principal, i.e. $I=(p)$ for some polynomial $p$. Write $\tilde{p}(w)=p(a,w)$ for $p\in\mathbb{C}[z,w]$, and $\tilde{f}(w)=f(a,w)$ for $f\in\mathcal{A}_\alpha(\mathbb{D}^2)$. Let $M_{\tilde{p}}$ denote the multiplication operator on $\mathcal{A}_\alpha(\mathbb{D})$ with symbol $\tilde{p}$.

    Since $Z(I)\cap(\{a\}\times\mathbb{T})=\emptyset$, it holds $\tilde{p}(w)\neq0,~\forall w\in\mathbb{T}$. List the roots of $\tilde{p}$ in $\mathbb{D}$ as $\{w_1,\ldots,w_n\}$, multiplicity being counted. Set $B(w)=\prod_{i=1}^n\frac{w_i-w}{1-\bar{w}_iw}$, then $B$ is a Blaschke product such that $\frac{\tilde{p}}{B}$ is rational and invertible on $\mathbb{T}$.

    For $f\in\mathcal{A}_\alpha(\mathbb{D}^2)$ we set
    $$Af=f-M_pM_{\frac{B}{\tilde{p}}(w)}M_{B(w)}^*\tilde{f}.$$
    For $f\in I=(p)$, suppose $f=pq$ then $\tilde{f}=\tilde{p}\tilde{q}$. It follows that
    \begin{eqnarray*}
        Af(z,w)&=&f(z,w)-p(z,w)\cdot\frac{B}{\tilde{p}}(w)\cdot\frac{\tilde{p}}{B}(w)\cdot\tilde{q}(w)\\
        &=&f(z,w)-p(z,w)\tilde{q}(w)\\
        &=&p(z,w)[q(z,w)-q(a,w)].
    \end{eqnarray*}
    Consequently $Af\in I$. For $w\in\mathbb{C}$, it holds
    $$Af(a,w)=p(a,w)[q(a,w)-q(a,w)]=0,$$
    which derives $Af\in I\cap(z-a)$. By $$I\cap(z-a)=(z-a)I\subset[(z-a)I]_\alpha,$$
    we get $Af\in[(z-a)I]_\alpha$. Consequently $A[I]_\alpha\subset[(z-a)I]_\alpha$. Suppose there exists $g\in[I]_\alpha^\perp$ such that $S_{z-a}g=0$. It means $(z-a)g\in[I]_\alpha$, then by definition
    $$A\big((z-a)g\big)=(z-a)g\in[(z-a)I]_\alpha.$$ Since $M_{z-a}$ on $\mathcal{A}_\alpha(\mathbb{D}^2)$ is bounded below, $M_{z-a}[I]_\alpha$ is closed and hence $M_{z-a}[I]_\alpha=[(z-a)I]_\alpha$. Then the injectivity of $M_{z-a}$ ensures $g\in[I]_\alpha$, which proves the injectivity of $S_{z-a}$.

    In the general case, let $I=pL$ be the Beurling decomposition, where $p$ is the greatest common divisor of all the polynomials in $I$, which is unique in the sense of a scalar multiplier, and $L$ is a finite codimensional ideal. If $g\in\ker S_{z-a}$, then $(z-a)g\in[I]_\alpha\subset[p]_\alpha$. By the proven injectivity of $S_{z-a}^{[p]_\alpha^\perp}$ we have $g\in[p]_\alpha$, and therefore from
    $$\ker S_{z-a}\subset[I]_\alpha^\perp\cap[p]_\alpha=[p]_\alpha\ominus[I]_\alpha$$
    we get
    $$\dim\ker S_{z-a}\leq\dim~(p)/I=\dim~\mathbb{C}[z,w]/L<\infty.$$
\end{proof}
The following lemma shows, roughly speaking that, $[I]_\alpha$ is essentially primary.
\begin{lem}\label{lem:injectivity2}
    Let $I$ be an ideal and $p\in\mathbb{C}[z,w]$, such that each component of $Z(p)$ intersects $\mathbb{D}^2$, and $Z(p)\cap Z(I)\cap\partial\mathbb{D}^2=\emptyset$, then $\ker S_p$ is finite dimensional. If in addition, $I$ is principal, then $S_p$ is injective.
\end{lem}
\begin{proof}
    By assumption, $Z(p)$ and $Z(I)$ do not have common algebraic components, and therefore $Z(p)\cap Z(I)$ is a finite set. Suppose
    $$Z(p)\cap Z(I)=\{(z_i,w_i):i=1,\ldots,n\},$$
    then by $Z(p)\cap Z(I)\cap\partial\mathbb{D}^2=\emptyset$, for each $i$ we have  $|z_i|\neq1$ and $|w_i|\neq1$. From Hilbert's Nullstellensatz, there is a nonconstant polynomial $q(z,w)=\prod_{i=1}^{n}(z-z_i)^{n_i}(w-w_i)^{m_i}$ in the ideal $(p)+I$, where $n_i,m_i$ are non-negative integers, such that $n_i=0$ whenever $|z_i|=1$, and $m_i=0$ whenever $|w_i|=1$.

    Let $I=p_0L$ be the Beurling decomposition, then by Lemma \ref{lem:injectivity1}, $S_q^{[p_0]_\alpha^\perp}$ is injective. From $q\in(p)+I$, there exist $f\in\mathbb{C}[z,w]$ and $q_1\in I$ such that $q=fp+q_1$. Then from $$S_q^{[p_0]_\alpha^\perp}=S_{f}^{[p_0]_\alpha^\perp}\cdot S_p^{[p_0]_\alpha^\perp}$$
    we get the injectivity of $S_p^{[p_0]_\alpha^\perp}$. When $I$ is principal, the injectivity of $S_p=S_p^{[p_0]_\alpha^\perp}$ has already been proved. For a general ideal $I$, if $g\in\ker S_p$ then $pg\in[I]_\alpha$, which belongs to $[p_0]_\alpha$. The injectivity of $S_p^{[p_0]_\alpha^\perp}$ gives $g\in[p_0]_\alpha$, and consequently
\begin{equation}\label{eq:injectivity}
    \ker S_p\subset[I]_\alpha^\perp\cap[p_0]_\alpha=[p_0]_\alpha\ominus[I]_\alpha
\end{equation}
is of finite dimension.
\end{proof}
For an operator $A\in B(H)$, $A$ is called an essential isometry if $\Id_H-A^*A$ is compact. Symmetrically, $A$ is essentially co-isometric if $A^*$ is essentially isometric. Moreover if $\Id_H-A^*A\in\mathfrak{L}^{(1,\infty)}$, $A$ is called $(1,\infty)$-essentially isometric; and if $A^*$ is $(1,\infty)$-essentially isometric, then $A$ is $(1,\infty)$-essentially co-isometric. If both $\Id_H-A^*A$ and $\Id_H-AA^*$ belong to $\mathfrak{L}^{(1,\infty)}$, then $A$ is called $(1,\infty)$-essentially unitary.

With these terminologies, we have the following result.
\begin{prop}\label{prop:enormal}
   Let $I\subset\mathbb{C}[z,w]$ be a distinguished ideal, and $\alpha\geq-1$. If $S_w$ is $(1,\infty)$-essentially co-isometric, then it is $(1,\infty)$-essentially unitary. Consequently, $P_{[I]_\alpha}M_wP_{[I]_\alpha^\perp}\in\mathfrak{L}^{(2,\infty)}$ and $[S_w^*,S_w]\in\mathfrak{L}^{(1,\infty)}$.
\end{prop}
\begin{proof}
    By assumption $S_wS_w^*$ is Fredholm. Consequently $\ker S_w^*=\ker S_wS_w^*$ is finite dimensional, and $\mathrm{ran}~S_w^*$ is closed. By Lemma \ref{lem:injectivity1}, $\ker S_w$ is finite dimensional, and therefore ensures the Fredholmness of $S_w^*$.

    Let $S_w=W|S_w|$ be the polar decomposition of $S_w$, where $|S_w|=(S_w^*S_w)^{1/2}$ and $W$ is the partial isometry that maps $(\ker S_w)^\perp$ unitarily onto $\mathrm{ran}~S_w$. Then we have
    $$\Id_{[I]_\alpha^\perp}-S_w^*S_w=W^*(\Id_{[I]_\alpha^\perp}-S_wS_w^*)W+(\Id_{[I]_\alpha^\perp}-W^*W)\in\mathcal{L}^{(1,\infty)},$$
    which proves that $S_w$ is essentially unitary.

    The remaining conclusions of the proposition follow from
    $$P_{[I]_\alpha^\perp}M_w^*P_{[I]_\alpha}M_wP_{[I]_\alpha^\perp}=P_{[I]_\alpha^\perp}(\Id_{[I]_\alpha^\perp}-S_w^*S_w)P_{[I]_\alpha^\perp}$$
    and
    $$[S_w^*,S_w]=(\Id_{[I]_\alpha^\perp}-S_wS_w^*)-(\Id_{[I]_\alpha^\perp}-S_w^*S_w).$$
\end{proof}
To prove the essential normality of $[I]_\alpha^\perp$, by Proposition \ref{prop:enormal} it suffices to prove both $S_z$ and $S_w$ are essentially co-isometric, and this will be the aim of the following two sections.

\section{Grassmannian quotient modules}\label{sec:Grassmannian}
~~~~Let $I\subset\mathbb{C}[z,w]$ be a distinguished ideal. Recall that for the distinguished variety $V=Z(I)\cap\mathbb{D}^2$,
$$m_V=\max_{\lambda\in\mathbb{D}}\card\{w\in\mathbb{D}:(\lambda,w)\in V\},~n_V=\max_{\lambda\in\mathbb{D}}\card\{z\in\mathbb{D}:(z,\lambda)\in V\}.$$
By \cite{AM1} there exists a rational $m_V\times m_V$-matrix inner function $\Psi(z)$ such that
$$V=\left\{(z,w)\in\mathbb{D}^2:\det\big(wI_{m_V}-\Psi(z)\big)=0\right\}.$$
To investigate the determinant $\det\big(wI_{m_V}-\Psi(z)\big)$, we introduce the construction of Grassmannian quotient module.

At first, we illustrate the main idea of the construction of Grassmaniann quotient module.  It is well-known that, the tensor product $\bigotimes\limits_{i=1}^{m_V}H^2(\mathbb{D})$ is unitarily equivalent to $H^2(\mathbb{D}^{m_V})$ via the mapping
$\bigotimes\limits_{i=1}^{m_V}f_i\mapsto\prod\limits_{i=1}^{m_V}f_i(z_i),$ and hence for convenience we identify $\bigotimes\limits_{i=1}^{m_V}H^2(\mathbb{D})$ with $H^2(\mathbb{D}^{m_V})$. Let $J_{m_V}$ be the ideal of $\mathbb{C}[z_1,\ldots,z_{m_V}]$ generated by $\{z_i-z_j:1\leq i,j\leq m_V\}$, then $[J_{m_V}]^\perp$ is unitarily equivalent to the weighted Bergman module $\mathcal{A}_{m_V-2}(\mathbb{D})$  \cite{Cla,FR1,FR2,WZ1} via
$$f(z_1,\ldots,z_{m_V})\mapsto f(z,\ldots,z),$$
 the restriction of $H^2(D^{m_V})$ on the diagonal $Z(J_{m_V})$.

Consequently, for $p_i\in\mathbb{C}[z]$ the action of $\bigotimes\limits_{i=1}^{m_V}p_i$  on $[J_{m_V}]^\perp$ is unitarily equivalent to the multiplication by $\prod\limits_{i=1}^{m_V}p_i(z)$ on $\mathcal{A}_{m_V-2}(\mathbb{D})$. Therefore $\mathcal{A}_{m_V-2}(\mathbb{D})$ is exactly the $(m_V-1)$-th tensor product of $H^2(\mathbb{D})$ on the ring $\mathbb{C}[z]$, and readers are referred to \cite{At} for tensor product of modules.  Following this idea, we construct the tensor product of general vector-valued Hilbert modules, and then define the Grassmannian quotient module as the maximal anti-symmetric quotient module of this tensor product.

Let $H$ be a Hilbert space, $\mathcal{A}$ be a commutative unital algebra, and $\psi:\mathcal{A}\to B(H)$ be an algebraic homomorphism. Then $\psi$ defines a module action on $H$, which makes $H$ into a Hilbert $\mathcal{A}$-module \cite{DP}. For simplicity, we write
$$ah=\psi(a)h,~a\in\mathcal{A},~h\in H.$$

Fix an integer $m>0$, write $H^m=\bigoplus\limits_{i=1}^mH_i$, where each $H_i$ is a copy of $H$. We identify the Hilbert space
$\tensorspace$ with $\bigoplus\limits_{\gamma\in\{1,\ldots,m\}^m}\bigotimes\limits_{i=1}^mH_{\gamma_i}$ via the mapping
$$\bigotimes_{i=1}^m\left(
                      \begin{array}{c}
                        h_{i,1} \\
                        \vdots \\
                        h_{i,m} \\
                      \end{array}
                    \right)\mapsto\bigoplus_{\gamma\in\{1,\ldots,m\}^m}\bigotimes_{i=1}^mh_{i,\gamma_i}.$$
                     There is a nature $\big(\bigotimes\limits_{i=1}^m\mathcal{A}\big)$-module structure on $\tensorspace$ defined by
$$\left(\bigotimes_{i=1}^ma_i\right)\left(\bigotimes_{i=1}^mh_i\right)=\bigotimes_{i=1}^m\big(a_ih_i\big),~a_i\in\mathcal{A},~h_i\in H_i^m.$$
For $f\in\tensorspace$ and $\gamma\in\{1,\ldots,m\}^m$, we denote by $f_\gamma$ the component of $f$ in $\bigotimes\limits_{i=1}^mH_{\gamma_i}$.

For matrix $A\in M_m(\mathcal{A})$ and $i=1,\ldots,m$, let
$$M_A:H^m\to H^m,f\mapsto Af,$$
and set
$$M_{A,i}=\Id_{H_1^m}\otimes\cdots\otimes\Id_{H_{i-1}^m}\otimes M_A\otimes\Id_{H_{i+1}^m}\otimes\cdots\otimes\Id_{H_m^m}:\tensorspace\to\tensorspace.$$
For $a\in\mathcal{A}$, we abbreviate $M_{aI_m,i}$ as $M_{a,i}$. The following three lemmas are obvious.
\begin{lem}
    For $A\in M_m(\mathcal{A})$ and $1\leq i\leq m$, it holds
    $$M_{A,i}^*=\Id_{H_1^m}\otimes\cdots\otimes\Id_{H_{i-1}^m}\otimes M_A^*\otimes\Id_{H_{i+1}^m}\otimes\cdots\otimes\Id_{H_m^m}.$$
\end{lem}
\begin{lem}\label{lem:commutativity}
    For $A,B\in M_m(\mathcal{A})$ and $i\neq j$, it holds
    $$[M_{A,i},M_{B,j}]=[M_{A,i}^*,M_{B,j}]=0,$$
    that is, $M_{A,i}$ and $M_{B,j}$ are doubly commuting.
\end{lem}
An isometry $S$ is called pure if $(SOT)\lim_{n\to\infty}S^{*n}=0$.
\begin{lem}\label{lem:isometrytensor}
    If $A\in M_m(\mathcal{A})$ is an (pure) isometry, $1\leq i\leq m$, then $M_{A,i}$ is also an (pure) isometry.
\end{lem}
For $A\in M_m(\mathcal{A})$, set $T_A=\prod_{i=1}^mM_{A,i}$. It is routine to verify
\begin{equation}\label{eq:TA}
    (T_Af)_\gamma=\sum_{\beta\in\{1,\ldots,m\}^m}\big(\bigotimes_{i=1}^ma_{\gamma_i,\beta_i}\big)f_\beta,\quad\gamma\in\{1,\ldots,m\}^m.
\end{equation}
Repeating application of Lemma \ref{lem:isometrytensor} gives:
\begin{cor}\label{cor:isometrytensor}
    If $A\in M_m(\mathcal{A})$ is isometric, then $T_A:\tensorspace\to\tensorspace$ is also isometric.
\end{cor}
For $\gamma=(\gamma_1,\ldots,\gamma_m)\in\{1,\ldots,m\}^m$ and $\sigma\in\mathcal{S}_m$, the symmetric group of degree $m$, denote by $\sigma(\gamma)=(\gamma_{\sigma(1)},\ldots,\gamma_{\sigma(m)})$. For $f\in\tensorspace$, define $\sigma(f)\in\tensorspace$ by
$$\sigma(f)_\gamma=f_{\sigma(\gamma)},~\gamma\in\{1,\ldots,m\}^m.$$
For simplicity, we write $f_{\sigma(1,\ldots,m)}$ as $f_\sigma$ in short.

On the tensor product $\tensorspace$, it holds $M_{c,i}=M_{c,j}$ for scalar $c$ and $1\leq i,j\leq m$. To represent the determinant $\det(wI_m-\Psi(z))$, we need to find a maximal quotient module $\mathcal{N}$ of $\tensorspace$ satisfying the following two conditions:
\begin{itemize}
    \item[(a)] $S_{a,i}^\mathcal{N}=S_{a,j}^\mathcal{N}$ for any $a\in\mathcal{A}$, where $S_{a,i}^\mathcal{N}=P_\mathcal{N}M_{a,i}\mid_\mathcal{N}$;
    \item[(b)] $P_\mathcal{N}\big(\bigotimes\limits_{i=1}^mf_i\big)=0$ whenever $f_1,\ldots,f_m\in H^m$ and there exist $i\neq j$ such that $f_i=f_j$.
\end{itemize}
Condition (a) requires
\begin{itemize}
    \item[(a')] $(M_{a,i}-M_{a,j})f\perp\mathcal{N},~\forall f\in\tensorspace,~a\in\mathcal{A},~1\leq i,j\leq m,$
\end{itemize}
and condition (b) is equivalent to
\begin{itemize}
    \item[(b')] for $f_1,\ldots,f_m\in H^m$ and $i<j$, it holds
    $$P_\mathcal{N}(f_1\otimes\cdots\otimes f_m)=-P_\mathcal{N}(f_1\otimes\cdots\otimes f_{i-1}\otimes f_j\otimes f_{i+1}\otimes\cdots\otimes f_{j-1}\otimes f_i\otimes f_{j+1}\otimes\cdots\otimes f_m).$$
\end{itemize}
which is furthermore equivalent to
\begin{itemize}
    \item[(b")] $f-\sgn(\sigma)\sigma(f)\perp\mathcal{N},~\forall f\in\tensorspace,~\sigma\in\mathcal{S}_m.$
\end{itemize}
According to conditions (a') and (b"), let $\mathcal{M}_{Gr}(H,m)$ be the closed subspace of $\tensorspace$ spanned by
$$\left\{f-\sgn(\sigma)\sigma(f),~(M_{a,i}-M_{a,j})f:f\in\tensorspace,~\sigma\in\mathcal{S}_m,~a\in\mathcal{A},~1\leq i,j\leq m\right\}.$$
Observe that the module action $\bigotimes\limits_{i=1}^ma_i=\prod\limits_{i=1}^mM_{a_i,i}$ on $\tensorspace$ is component-wise, and hence commutes with $\sigma$. On the other hand, for $f\in\tensorspace$ we have
$$M_{a,j}\big(\bigotimes_{i=1}^ma_i\big)f=M_{a,j}\prod_{i=1}^mM_{a_i,i}f=\big(\prod_{i=1}^mM_{a_i,i}\big)M_{a,j}f=\big(\bigotimes_{i=1}^ma_i\big)M_{a,j}f,$$
and therefore $M_{a,j}$ commutes with the module actions of $\tensorspace$. Then it is routine to verify that $\mathcal{M}_{Gr}(H,m)$ is a submodule of $\tensorspace$. We call
$$\mathcal{N}_{Gr}(H,m)=\left(\tensorspace\right)\ominus\mathcal{M}_{Gr}(H,m)$$
the $m$-th Grassmannian quotient module associated with $H$. By definition, it clearly holds
\begin{equation}\label{eq:Atensor}
    S_{a,i}^{\mathcal{N}_{Gr}(H,m)}=S_{a,j}^{\mathcal{N}_{Gr}(H,m)},~a\in\mathcal{A},~1\leq i,j\leq m,
\end{equation}
and
\begin{equation}\label{eq:antisymmetricity}
    f=\sgn(\sigma)\sigma(f),~f\in\mathcal{N}_{Gr}(H,m),~\sigma\in\mathcal{S}_m
\end{equation}
as required.

For distinct $i,j\in\{1,\ldots,m\}$, denote by $\sigma_{i,j}\in\mathcal{S}_m$, the permutation that exchanges the $i$-th and $j$-th elements. If $\gamma\in\{1,\ldots,m\}^m$ such that $\gamma_i=\gamma_j$ for some $i\neq j$, then we have $\sigma_{i,j}(\gamma)=\gamma$. By (\ref{eq:antisymmetricity}), for $f\in\mathcal{N}_{Gr}(H,m)$ we have $f=-\sigma_{i,j}(f)$, which yields $f_\gamma=-f_\gamma$ and consequently $f_\gamma=0$. Therefore the only significant components of $f$ are those $f_\sigma,~\sigma\in\mathcal{S}_m$. Set $f_0=f_{(1,2,\ldots,m)}$, then $$f_\sigma=\sgn(\sigma)f_0,~\forall\sigma\in\mathcal{S}_m.$$

For the module $\bigotimes\limits_{i=1}^mH_i$ where each $H_i=H$, and $a\in\mathcal{A},~1\leq i\leq m$, we set
$$M_{a,i}=\Id_{H_1}\otimes\cdots\otimes\Id_{H_{i-1}}\otimes a\otimes\Id_{H_{i+1}}\otimes\cdots\otimes\Id_{H_m}.$$
Write $\mathcal{M}_0$ the submodule generated by
$$\{(M_{a,i}-M_{a,j})f:f\in\bigotimes_{i=1}^mH_i,a\in\mathcal{A},1\leq i,j\leq m\},$$
and $\mathcal{N}_0=\big(\bigotimes\limits_{i=1}^mH_i\big)\ominus\mathcal{M}_0$. From the definition of $\mathcal{N}_{Gr}(H,m)$, for $f\in\mathcal{N}_{Gr}(H,m)$ we have $f_0\in\mathcal{N}_0$. Conversely, for each $f_0\in\mathcal{N}_0$, the element $f\in\tensorspace$ defined by $f_\sigma=\sgn(\sigma)f_0,~\sigma\in\mathcal{S}_m$ belongs to $\mathcal{N}_{Gr}(H,m)$. Consequently, $\mathcal{N}_{Gr}(H,m)$ is unitarily isomorphic to $\mathcal{N}_0$, via the mapping
$$\Theta:\mathcal{N}_{Gr}(H,m)\to\mathcal{N}_0,f\mapsto\sqrt{m!}f_0,$$
and hence
\begin{equation}\label{eq:piequivalence}
    M_{a,i}^*\Theta=\Theta M_{a,i}^*,~\forall a\in\mathcal{A},~1\leq i\leq m.
\end{equation}
\begin{lem}\label{lem:permute}
    Given $\sigma\in\mathcal{S}_m$, let $A$ be the $m\times m$-matrix such that $Ax=(x_{\sigma(1)},\ldots,x_{\sigma(m)})^T$ for all vectors $x\in\mathbb{C}^m$, then
    $$(T_Af)_\gamma=f_{(\sigma(\gamma_1),\ldots,~\sigma(\gamma_m))}$$
    for all $f\in\tensorspace$ and $\gamma\in\{1,\ldots,m\}^m$.
\end{lem}
\begin{proof}
    Clearly $A$ is the matrix defined by $a_{i,j}=\delta_{\sigma(i)}(j)$, where $\delta$ is the Kronecker symbol. Then by (\ref{eq:TA}) we have
    \begin{eqnarray*}
        (T_Af)_\gamma&=&\sum_{\beta\in\{1,\ldots,m\}^m}\big(\prod_{i=1}^ma_{\gamma_i,\gamma_i}\big)f_\beta\\
        &=&\prod_{i=1}^ma_{\gamma_i,\sigma(\gamma_i)}f_{(\sigma(\gamma_1),\ldots,\sigma(\gamma_m))}\\
        &=&f_{(\sigma(\gamma_1),\ldots,\sigma(\gamma_m))}.
    \end{eqnarray*}
\end{proof}
Denote by
$$S_A^{\mathcal{N}_{Gr}(H,m)}=P_{\mathcal{N}_{Gr}(H,m)}T_A\mid_{\mathcal{N}_{Gr}(H,m)},~A\in M_m(\mathcal{A}),$$
and
$$S_{a,i}^{\mathcal{N}_0}=P_{\mathcal{N}_0}M_{a,i}\mid_{\mathcal{N}_0},~a\in\mathcal{A},~1\leq i\leq m.$$
By the definition of $\mathcal{N}_0$, for $a\in\mathcal{A}$ and any $i$ and $j$ we have
\begin{equation}\label{eq:MaiMaj}
    M_{a,i}^*\mid_{\mathcal{N}_0}=M_{a,j}^*\mid_{\mathcal{N}_0},
\end{equation}
which in turn gives $S_{a,i}^{\mathcal{N}_0}=S_{a,j}^{\mathcal{N}_0}$. We write $S_a^{\mathcal{N}_0}$ for $S_{a,i}^{\mathcal{N}_0},~a\in\mathcal{A},~1\leq i\leq m$. Then the mapping $a\mapsto S_a^{\mathcal{N}_0}$ defines an $\mathcal{A}$-module structure on  $\mathcal{N}_0$.

Grassmiannian quotient module enjoys the following natural and intrinsic property.
\begin{thm}\label{thm:Grassmannian}
    $\mathcal{N}_{Gr}(H,m)$ is invariant for $T_A^*,~\forall A\in M_m(\mathcal{A})$, and $\Theta T_A^*~f=\big(S_{\det A}^{\mathcal{N}_0}\big)^*\Theta~f$ for $f\in\mathcal{N}_{Gr}(H,m)$.
\end{thm}
\begin{proof}
    If $f\in\mathcal{N}_{Gr}(H,m)$, write $f_0=f_{(1,\ldots,m)}\in\mathcal{N}_0$, then $f_\gamma=0$ when $\gamma$ is not a permutation of $(1,2,\ldots,m)$, and $f_\sigma=\sgn(\sigma) f_0$ if $\sigma\in\mathcal{S}_m$. By the definition of $\mathcal{N}_0$, for $f\in\mathcal{N}_{Gr}(H,m)$ we have
    $$\big(S_a^{\mathcal{N}_0}\big)^*f_\sigma=M_{a,i}^*f_\sigma=M_{a,j}^*f_\sigma,~\forall a\in\mathcal{A},~\sigma\in\mathcal{S}_m,~1\leq i,j\leq m.$$
    By Lemma \ref{lem:commutativity} it obviously holds
    \begin{equation}\label{eq:lem:determinant1}
        (M_{a,i}^*-M_{a,j}^*)T_A^*f=(M_{a,i}^*-M_{a,j}^*)\prod_{k=1}^mM_{A,k}^*f=\prod_{k=1}^mM_{A,k}^*(M_{a,i}^*-M_{a,j}^*)f=0.
    \end{equation}
    For $f\in\mathcal{N}_{Gr}(H,m),~g\in\tensorspace$ and $\sigma\in\mathcal{S}_m$, by (\ref{eq:TA}) and (\ref{eq:MaiMaj}) we have
    \begin{eqnarray*}
        \langle f,T_A\sigma(g)\rangle&=&\sum_{\rho\in\mathcal{S}_m}\left\langle f_\rho,\big(T_A\sigma(g)\big)_\rho\right\rangle\\
        &=&\sum_{\rho\in\mathcal{S}_m}\left\langle f_\rho,\sum_{\beta\in\{1,\ldots,m\}^m}\big(\bigotimes_{j=1}^ma_{\rho(j),\beta_j}\big)\sigma(g)_\beta\right\rangle\\
        &=&\sum_{\rho\in\mathcal{S}_m}\sum_{\beta\in\{1,\ldots,m\}^m}\left\langle f_\rho,\prod_{j=1}^mM_{a_{\rho(j),\beta_j},j}g_{\sigma(\beta)}\right\rangle\\
        &=&\sum_{\rho\in\mathcal{S}_m}\sum_{\beta\in\{1,\ldots,m\}^m}\left\langle f_\rho,\prod_{j=1}^mM_{a_{\rho\sigma(j),\sigma(\beta)_j},j}g_{\sigma(\beta)}\right\rangle\\
        &=&\sum_{\rho\in\mathcal{S}_m}\left\langle\sgn(\sigma)f_{\rho\sigma},(T_Ag)_{\rho\sigma}\right\rangle\\
        &=&\left\langle f,\sgn(\sigma)T_Ag\right\rangle,
    \end{eqnarray*}
    which induces
    \begin{eqnarray}\label{eq:lem:determinant2}
        \langle T_A^*f,g-\sgn(\sigma)\sigma(g)\rangle&=&\langle f,T_Ag-\sgn(\sigma)T_A\sigma(g)\rangle\notag\\
        &=&\langle f,T_Ag-T_Ag\rangle\notag\\
        &=&0.
    \end{eqnarray}
    From (\ref{eq:lem:determinant1}) and (\ref{eq:lem:determinant2}) we conclude $T_A^*f\in\mathcal{N}_{Gr}(H,m)$, i.e. $\mathcal{N}_{Gr}(H,m)$ is invariant for $T_A^*$.

    For $f,g\in\mathcal{N}_{Gr}(H,m)$, by (\ref{eq:piequivalence}) we have
    \begin{eqnarray*}
        \left\langle\Theta T_A^*f,\Theta g\right\rangle&=&\left\langle f,T_Ag\right\rangle\\
        &=&\sum_{\sigma\in\mathcal{S}_m}\left\langle f_\sigma,\sum_{\rho\in\mathcal{S}_m}\big(\bigotimes_{j=1}^ma_{\sigma(j),\rho(j)}\big)g_\rho\right\rangle\\
        &=&\sum_{\sigma\in\mathcal{S}_m}\left\langle f_0,\sum_{\rho\in\mathcal{S}_m}\sgn(\sigma\rho)\prod_{j=1}^mM_{a_{\sigma(j),\rho(j)},j}g_0\right\rangle\\
        &=&\sum_{\sigma\in\mathcal{S}_m}\left\langle f_0,\sum_{\rho\in\mathcal{S}_m}\sgn(\sigma\rho)\prod_{j=1}^mS_{a_{\sigma(j),\rho(j)}}^{\mathcal{N}_0}g_0\right\rangle\\
        &=&\frac{1}{m!}\sum_{\sigma\in\mathcal{S}_m}\left\langle\sum_{\rho\in\mathcal{S}_m}\sgn(\rho\sigma^{-1})\prod_{j=1}^m\big(S_{a_{j,\rho\sigma^{-1}(j)}}^{\mathcal{N}_0}\big)^*\Theta f,\Theta g\right\rangle\\
        &=&\frac{1}{m!}\sum_{\sigma\in\mathcal{S}_m}\left\langle\big(S_{\det A}^{\mathcal{N}_0}\big)^*\Theta f,\Theta g\right\rangle\\
        &=&\left\langle\big(S_{\det A}^{\mathcal{N}_0}\big)^*\Theta f,\Theta g\right\rangle,
    \end{eqnarray*}
    which completes the proof of the theorem.
\end{proof}

\section{Essential normality of quotient module $[I]_\alpha^\perp$}
~~~~Now we turn attention back to the modules over the bidisk. Set $\mathcal{A}=\mathbb{C}[z,w]$ and $H=H^2(\mathbb{D}^2)$. By the preceding section, we define the Grassmannian quotient module $\mathcal{N}_{Gr}(H^2(\mathbb{D}^2),m)$. In this case, $\mathcal{M}_0$ is the submodule of $\bigotimes\limits_{i=1}^mH^2(\mathbb{D}^2)=H^2(\mathbb{D}^{2m})$ generated by
$$\{z_i-z_j,~w_i-w_j:1\leq i,j\leq m\}.$$
Let $J_m$ be the ideal of $\mathbb{C}[z_1,\ldots,z_m]$ generated by $\{z_i-z_j:1\leq i,j\leq m\}$, and $[J_m]^\perp$ the associated quotient module in $H^2(\mathbb{D}^m)$. As mentioned in the second paragraph of Section \ref{sec:Grassmannian}, $[J_m]^\perp$ is isometrically isomorphic to $\mathcal{A}_{m-2}^2(\mathbb{D})$ via $f(z_1,\ldots,z_m)\mapsto f(z,\ldots,z)$. By the definition of $\mathcal{N}_0$ we have $\mathcal{N}_0=[J_m]^\perp\otimes[J_m]^\perp$, which is isomorphic to $\mathcal{A}_{m-2}^2(\mathbb{D}^2)=\mathcal{A}_{m-2}^2(\mathbb{D})\otimes\mathcal{A}_{m-2}^2(\mathbb{D})$. Therefore we may redefine
\begin{eqnarray*}
    \Theta:\mathcal{N}_{Gr}(H^2(\mathbb{D}^2),m)&\to&\mathcal{A}_{m-2}^2(\mathbb{D}^2)\\
    f(z_1,\ldots,z_m,w_1,\ldots,w_m)&\mapsto&\sqrt{m!}\sgn(\sigma)f_\sigma(z,\ldots,z,w,\ldots,w).
\end{eqnarray*}
\begin{lem}\label{lem:fg}
    For a rational inner $\Psi(z)\in M_m(H^\infty(\mathbb{D}))$, the subspace
    $$\bigcap_{i=1}^m\big(\ker M_{\Psi(z),i}^*\bigcap\ker M_{wI_m,i}^*\big)$$
    is of finite dimension.
\end{lem}
\begin{proof}
    By Lemma \ref{lem:commutativity} and Lemma \ref{lem:isometrytensor}, $\{M_{\Psi(z),i},M_{wI_m,i}:i=1,\ldots,m\}$ is a doubly commuting family of isometries. For each $i$, $\ker M_{\Psi(z)}^*\cap\ker M_{wI_m}^*$ consists of the vector-valued functions $f$ in $\overline{\mathrm{span}}~\{z^n:n=0,1,\ldots\}\otimes\mathbb{C}^m$ such that $M_{\Psi(z)}^*f=0$. Since $\Psi(z)$ is rational, $\ker M_{\Psi(z)}^*\cap\ker M_{wI_m}^*$ is finite dimensional. Therefore
    $$\bigcap_{i=1}^m\big(\ker M_{\Psi(z),i}^*\bigcap\ker M_{wI_m,i}^*\big)=\bigotimes_{i=1}^m\big(\ker M_{\Psi(z)}^*\bigcap\ker M_{wI_m}^*\big)$$
    is finite dimensional.
\end{proof}
The following proposition is the key point of this paper.
\begin{prop}\label{prop:main}
     Let $I\subset\mathbb{C}[z,w]$ be a distinguished ideal, and suppose that $\Psi(z)\in M_m(H^\infty(\mathbb{D}))$ is pure rational inner, such that $$Z(I)\subset Z\big(\det\big(wI_m-\Psi(z)\big)\big).$$
     Then on the quotient weighted Bergman module $[I]_{m-2}^\perp$, $S_w$ is $(1,\infty)$-essentially unitary, and $[S_w^*,S_w]\in\mathfrak{L}^{(1,\infty)}$.
\end{prop}
\begin{proof}
    Let $p(z,w)$ be the numerator of the rational function $\det\big(wI_m-\Psi(z)\big)$, then by Hilbert's Nullstellensatz, there is a natural number $N$ such that $p(z,w)^N\in I$, which gives
    $$\det\big(wI_m-\Psi(z)\big)^N\in[I]_{m-2}.$$
    Set
    $$K=\bigcap_{i=1}^m\big(\ker M_{\Psi(z),i}^*\bigcap\ker M_{wI_m,i}^*\big),$$
    then by Lemma \ref{lem:fg} we have $\dim K<\infty$. By Lemma \ref{lem:commutativity} and Lemma \ref{lem:isometrytensor}, $$\{M_{\Psi(z),i},M_{wI_m,i}:i=1,\ldots,m\}$$
    is a doubly commuting family of pure isometries. Applying the Wold decomposition we have
    $$\tensorspace=\bigoplus_{\alpha,\beta\in\mathbb{N}^m}\prod_{i=1}^mM_{wI_m,i}^{\alpha_i}M_{\Psi(z),i}^{\beta_i}K.$$
    For
    $$\sum_{\alpha,\beta\in\mathbb{N}^m}\prod_{i=1}^mM_{wI_m,i}^{\alpha_i}M_{\Psi(z),i}^{\beta_i}f_{\alpha,\beta}\in\tensorspace,$$
    where each $f_{\alpha,\beta}\in K$ and $\sum_{\alpha,\beta\in\mathbb{N}^m}\|f_{\alpha,\beta}\|^2<\infty,$ define
    $$U\big(\sum_{\alpha,\beta\in\mathbb{N}^m}\prod_{i=1}^mM_{wI_m,i}^{\alpha_i}M_{\Psi(z),i}^{\beta_i}f_{\alpha,\beta}\big)=\sum_{\alpha,\beta\in\mathbb{N}^m}\big(\prod_{i=1}^mw_i^{\alpha_i}z_i^{\beta_i}\big)\otimes f_{\alpha,\beta},$$
    then $U:\tensorspace\to H^2(\mathbb{D}^{2m})\otimes K$ is a unitary operator such that
    \begin{eqnarray}\label{eq:unitary}
        M_{wI_m,i}&=&U^*(M_{w_i}\otimes\Id_K)U,\label{eq:unitary1}\\
        M_{\Psi(z),i}&=&U^*(M_{z_i}\otimes\Id_K)U.\notag
    \end{eqnarray}
    Let $\mathcal{M}_1\subset H^2(\mathbb{D}^{2m})$ be the homogeneous submodule generated by
    $$\{w_i-w_j,z_i-z_j,\prod_{k=1}^m(w_k-z_k)^N:1\leq i,j\leq m\},$$
    and $\mathcal{N}_1=\mathcal{M}_1^\perp$ the associated quotient module. Recall that for the ideal $J_m\subset\mathbb{C}[z_1,\ldots,z_m]$ generated by $\{z_i-z_j:1\leq i,j\leq m\}$, the quotient Hardy module $[J_m]^\perp\subset H^2(\mathbb{D}^m)$ is isometrically isomorphic to $\mathcal{A}_{m-2}^2(\mathbb{D})$ via $f(z_1,\ldots,z_m)\mapsto f(z,\ldots,z)$. Therefore the unitary isomorphism defined by
    \begin{eqnarray*}
        [J_m]^\perp\otimes[J_m]^\perp&\to&\mathcal{A}_{m-2}(\mathbb{D}^2),\\
        f(z_1,\ldots,z_m)\otimes g(w_1,\ldots,w_d)&\mapsto&f(z,\ldots,z)g(w,\ldots,w)
    \end{eqnarray*}
    maps $\mathcal{N}_1$ onto the quotient module $[(z-w)^{mN}]_{m-2}^\perp\subset\mathcal{A}_{m-2}^2(\mathbb{D})$. The proof of \cite[Corollary 2.7]{WZ1} for quotient Hardy module is also valid to quotient modules of $\mathcal{A}_\alpha(\mathbb{D}^d)(\alpha>-1)$, and hence we have
    $$\Id_{[(z-w)^{mN}]_{m-2}^\perp}-\big(S_w^{[(z-w)^{mN}]_{m-2}^\perp}\big)\big(S_w^{[(z-w)^{mN}]_{m-2}^\perp}\big)^*\in\mathfrak{L}^{(1,\infty)},$$
    which in turn gives
    \begin{equation}\label{eq:eisometry1}
        \Id_{\mathcal{N}_1}-\big(S_{w_i}^{\mathcal{N}_1}\big)\big(S_{w_i}^{\mathcal{N}_1}\big)^*\in\mathfrak{L}^{(1,\infty)},1\leq i\leq m.
    \end{equation}
    Set $\mathcal{M}_2=U^*(\mathcal{M}_1\otimes K)$ and $\mathcal{N}_2=U^*(\mathcal{N}_1\otimes K)$, then since $U$ is unitary, by (\ref{eq:unitary1}) and (\ref{eq:eisometry1}) we obtain
    \begin{equation}\label{eq:eisometry2}
        \Id_{\mathcal{N}_2}-\big(S_{wI_m,i}^{\mathcal{N}_2}\big)\big(S_{wI_m,i}^{\mathcal{N}_2}\big)^*\in\mathfrak{L}^{(1,\infty)}.
    \end{equation}
    Write $\mathcal{M}_3$ the submodule of $\tensorspace$ generated by
    $$\{(M_{wI_m,i}-M_{wI_m,j})f,(M_{\Psi(z),i}-M_{\Psi(z),j})f,T_{wI_m-\Psi(z)}^Nf:f\in\tensorspace,1\leq i,j\leq m\}$$
    and $\mathcal{N}_3=\mathcal{M}_3^\perp$. It is clear that $\mathcal{M}_2\subset\mathcal{M}_3$, which implies $\mathcal{N}_3\subset\mathcal{N}_2$. On the other hand, from $\det\big(wI_m-\Psi(z)\big)^N\in[I]_{m-2}$ we get
    $$[I]_{m-2}^\perp\subset\ker M_{\det(wI_m-\Psi(z))}^{*N},$$
    and hence Theorem \ref{thm:Grassmannian} implies that
    $$\Theta^{-1}[I]_{m-2}^\perp\subset\left(\mathcal{N}_{Gr}(H,m)\cap\ker T_{wI_m-\Psi(z)}^{*N}\right)\subset\mathcal{N}_3\subset\mathcal{N}_2.$$
    Consequently (\ref{eq:eisometry2}) ensures that $$\Id_{\Theta^{-1}[I]_{m-2}^\perp}-\big(S_{wI_m,i}^{\Theta^{-1}[I]_{m-2}^\perp}\big)\big(S_{wI_m,i}^{\Theta^{-1}[I]_{m-2}^\perp}\big)^*\in\mathfrak{L}^{(1,\infty)},$$
    which together with (\ref{eq:piequivalence}) gives
    $$\Id_{[I]_{m-2}^\perp}-\big(S_w^{[I]_{m-2}^\perp}\big)\big(S_w^{[I]_{m-2}^\perp}\big)^*\in\mathfrak{L}^{(1,\infty)}.$$
    Now the desired conclusions follow from Proposition \ref{prop:enormal}.
\end{proof}
\begin{thm}\label{thm:main1}
    Let $I\subset\mathbb{C}[z,w]$ be a distinguished ideal, and $V$ be the distinguished variety $Z(I)\cap\mathbb{D}^2$. Then for each integer
    $$\alpha\geq\max\{m_V,n_V\}-2,$$
    the quotient module $[I]_\alpha^\perp$ is $(1,\infty)$-essentially normal. Moreover, both $S_z$ and $S_w$ on $[I]_\alpha^\perp$ are $(1,\infty)$-essentially unitary.
\end{thm}
\begin{proof}
    By \cite{AM1} there is a pure rational inner matrix $\Psi\in M_{m_V}(H^\infty(\mathbb{D}))$ such that $Z(I)$ is the zero set of $\det\big(wI_{m_V}-\Psi(z)\big)$. Set
    $$\Psi_1(z)=\left(
              \begin{array}{cc}
                \Psi(z) & O \\
                O & zI_{\alpha-m_V} \\
              \end{array}
            \right),$$
            then $\Psi_1$ is a pure rational inner $\alpha\times\alpha$-matrix, with
            $$Z(I)\subset Z\left(\det\big(wI_\alpha-\Psi_1(z)\big)\right).$$
            Similarly there is a pure rational inner $\alpha\times\alpha$-matrix $\Psi_2$, with
            $$Z(I)\subset Z\left(\det\big(zI_\alpha-\Psi_2(w)\big)\right).$$
            Then the desired conclusions follow from Proposition \ref{prop:main}.
\end{proof}
Actually, Theorem \ref{thm:main1} gives a necessary and sufficient condition for a variety to be distinguished.
\begin{thm}\label{thm:distinguishedvariety}
Let $V$ be an algebraic subvariety of $\mathbb D^2$. If there exists $\alpha\geq-1$, such that all the quotient modules $[I]_\alpha^\perp$ with $Z(I)\cap\mathbb{D}^2=V$ are essentially normal, then $V$ is a distinguished variety.
\end{thm}
\begin{proof}
    Suppose $V$ is not distinguished, then by definition there exists an accumulation point $(z_0,w_0)$ of $V$, which lies in $\partial(\mathbb{D}^2)\backslash\mathbb{T}^2$. Without loss of generality, assume $|z_0|<1$ and $|w_0|=1$.

    Let
    $$I_V=\{p\in\mathbb{C}[z,w]:p\mid_V=0\},$$
    the radical ideal associated to $V$, and take $I=I_V^2$. It is clear that $V=Z(I)\cap\mathbb{D}^2$. We shall deny the essential normality of $S_z$ on $[I]_\alpha^\perp$.

    For each $\lambda\in V$, it is obvious that the reproducing kernel
    $$K_\lambda(z,w)=\frac{1}{(1-\bar{\lambda}_1z)^{\alpha+2}(1-\bar{\lambda}_2w)^{\alpha+2}},~(z,w)\in\mathbb{D}^2$$
    belongs to $[I]_\alpha^\perp$.

    For $\mu\in\mathbb{D}^2$ we have
    \begin{eqnarray}\label{eq:kernelderivative1}
        \langle K_\mu,z\partial_zK_{\lambda}\rangle&=&\left\langle K_\mu,\frac{(\alpha+2)\bar{\lambda}_1z}{(1-\bar{\lambda}_1z)^{\alpha+3}(1-\bar{\lambda}_2w)^{\alpha+2}}\right\rangle\notag\\
        &=&\frac{(\alpha+2)\bar{\mu}_1\lambda_1}{(1-\bar{\mu}_1\lambda_1)^{\alpha+3}(1-\bar{\mu}_2\lambda_2)^{\alpha+2}}\notag\\
        &=&(z\partial_zK_\mu)(\lambda),
    \end{eqnarray}
    which indicates by linearity
    \begin{equation}\label{eq:kernelderivative2}
        \langle f,z\partial_zK_{\lambda}\rangle=(z\partial_zf)(\lambda),~\forall f\in\mathcal{A}_\alpha(\mathbb{D}^2).
    \end{equation}
    Then for each $f,g\in I_V$ we have
    $$\langle fg,z\partial_zK_{\lambda}\rangle=[z\partial_z(fg)](\lambda)=[z(f\partial_zg+g\partial_zf)](\lambda)=0,$$
    which gives $z\partial_zK_{\lambda}\perp I_V^2$, and consequently $z\partial_zK_{\lambda}\in[I]_\alpha^\perp$.

    By (\ref{eq:kernelderivative2}) we have
    \begin{eqnarray}\label{eq:kernelderivative3}
        \|z\partial_zK_{\lambda}\|^2&=&z\partial_z(z\partial_zK_{\lambda})(\lambda)\notag\\
        &=&(z\partial_zK_{\lambda}+z^2\partial_z^2K_\lambda)(\lambda)\notag\\
        &=&\frac{(\alpha+2)|\lambda_1|^2}{(1-|\lambda_1|^2)^{\alpha+3}(1-|\lambda_2|^2)^{\alpha+2}}+\frac{(\alpha+2)(\alpha+3)|\lambda_1|^4}{(1-|\lambda_1|^2)^{\alpha+4}(1-|\lambda_2|^2)^{\alpha+2}}\notag\\
        &=&\frac{(\alpha+2)|\lambda_1|^2[1+(\alpha+2)|\lambda_1|^2]}{(1-|\lambda_1|^2)^{\alpha+4}(1-|\lambda_2|^2)^{\alpha+2}}.
    \end{eqnarray}
    Let $k_\lambda=\frac{K_\lambda}{\|K_\lambda\|}$ be the normalized reproducing kernel. Then (\ref{eq:kernelderivative1}) derives
    \begin{eqnarray}\label{eq:kernelderivative4}
        |\langle z\partial_zK_{\lambda},k_\lambda\rangle|^2&=&\frac{1}{\|K_\lambda\|^2}|\langle z\partial_zK_{\lambda},K_\lambda\rangle|^2\notag\\
        &=&(1-|\lambda_1|^2)^{\alpha+2}(1-|\lambda_2|^2)^{\alpha+2}\cdot\left|\frac{(\alpha+2)|\lambda_1|^2}{(1-|\lambda_1|^2)^{\alpha+3}(1-|\lambda_2|^2)^{\alpha+2}}\right|^2\notag\\
        &=&\frac{(\alpha+2)^2|\lambda_1|^4}{(1-|\lambda_1|^2)^{\alpha+4}(1-|\lambda_2|^2)^{\alpha+2}}.
    \end{eqnarray}
    Set
    $$f_\lambda=\frac{z\partial_zK_{\lambda}-\langle z\partial_zK_{\lambda},k_\lambda\rangle k_\lambda}{\|z\partial_zK_{\lambda}-\langle z\partial_zK_{\lambda},k_\lambda\rangle k_\lambda\|},$$
    then in the case $\|z\partial_zK_{\lambda}-\langle z\partial_zK_{\lambda},k_\lambda\rangle k_\lambda\|$ is nonzero, $f_\lambda$ is a unit vector in $[I]_\alpha^\perp$, and is orthogonal to $K_\lambda$.

    It follows from (\ref{eq:kernelderivative3}) and (\ref{eq:kernelderivative4}) that
    \begin{eqnarray}\label{eq:kernelderivative5}
        \|z\partial_zK_{\lambda}-\langle z\partial_zK_{\lambda},k_\lambda\rangle k_\lambda\|^2&=&\|z\partial_zK_{\lambda}\|^2-|\langle z\partial_zK_{\lambda},k_\lambda\rangle|^2\notag\\
        &=&\frac{(\alpha+2)|\lambda_1|^2[1+(\alpha+2)|\lambda_1|^2]}{(1-|\lambda_1|^2)^{\alpha+4}(1-|\lambda_2|^2)^{\alpha+2}}-\frac{(\alpha+2)^2|\lambda_1|^4}{(1-|\lambda_1|^2)^{\alpha+4}(1-|\lambda_2|^2)^{\alpha+2}}\notag\\
        &=&\frac{(\alpha+2)|\lambda_1|^2}{(1-|\lambda_1|^2)^{\alpha+4}(1-|\lambda_2|^2)^{\alpha+2}}\notag\\
        &=&\frac{(\alpha+2)|\lambda_1|^2}{(1-|\lambda_1|^2)^2}\|K_\lambda\|^2,
    \end{eqnarray}
    and by (\ref{eq:kernelderivative1}) we find
    \begin{eqnarray}\label{eq:kernelderivative6}
        \langle zk_\lambda,z\partial_zK_{\lambda}-\langle z\partial_zK_{\lambda},k_\lambda\rangle k_\lambda\rangle&=&\langle zk_\lambda,z\partial_zK_{\lambda}\rangle-\langle z\partial_zK_{\lambda},k_\lambda\rangle\langle zk_\lambda,k_\lambda\rangle\notag\\
        &=&\frac{1}{\|K_\lambda\|}\big(\langle zK_\lambda,z\partial_zK_{\lambda}\rangle-\lambda_1\langle z\partial_zK_{\lambda},K_\lambda\rangle\big)\notag\\
        &=&\frac{1}{\|K_\lambda\|}\big(z\partial_z(zK_\lambda)-\lambda_1z\partial_zK_{\lambda}\big)(\lambda)\notag\\
        &=&\frac{1}{\|K_\lambda\|}\big(zK_\lambda+z^2\partial_zK_\lambda-\lambda_1z\partial_zK_{\lambda}\big)(\lambda)\notag\\
        &=&\frac{\lambda_1K_\lambda(\lambda)}{\|K_\lambda\|}\notag\\
        &=&\lambda_1\|K_\lambda\|.
    \end{eqnarray}
    In the case $\lambda_1\neq0$, (\ref{eq:kernelderivative5}) and (\ref{eq:kernelderivative6}) gives
    \begin{eqnarray*}
        \|S_zk_{\lambda}\|^2-\|S_z^*k_{\lambda}\|^2&\geq&|\langle S_zk_{\lambda},f_{\lambda}\rangle|^2+|\langle S_zk_{\lambda},k_{\lambda}\rangle|^2-\|S_z^*k_{\lambda}\|^2\\
        &=&\big|\big\langle zk_{\lambda},\frac{z\partial_zK_{\lambda}-\langle z\partial_zK_{\lambda},k_\lambda\rangle k_\lambda}{\|z\partial_zK_{\lambda}-\langle z\partial_zK_{\lambda},k_\lambda\rangle k_\lambda\|}\big\rangle\big|^2+|\lambda_1|^2-|\lambda_1|^2\\
        &=&\frac{|\langle zk_{\lambda},z\partial_zK_{\lambda}-\langle z\partial_zK_{\lambda},k_\lambda\rangle k_\lambda\rangle|^2}{\|z\partial_zK_{\lambda}-\langle z\partial_zK_{\lambda},k_\lambda\rangle k_\lambda\|^2}\\
        &=&\frac{|\lambda_1|^2\|K_\lambda\|^2}{\frac{(\alpha+2)|\lambda_1|^2}{(1-|\lambda_1|^2)^2}\|K_\lambda\|^2}\\
        &=&\frac{(1-|\lambda_1|^2)^2}{\alpha+2}.
    \end{eqnarray*}

    If there is a sequence $\{\lambda^{(n)}\}\subset V$ with $\lambda_1^{(n)}\neq0$ approaches $(z_0,w_0)$, then
    $$\liminf_{n\to\infty}(\|S_zk_{\lambda^{(n)}}\|^2-\|S_z^*k_{\lambda^{(n)}}\|^2)\geq\frac{(1-|z_0|^2)^2}{\alpha+2}>0,$$
    contracting to the compactness of $[S_z^*,S_z]$, since $k_{\lambda^{(n)}}$ converges weakly to $0$.

    If there does not exist such a sequence $\{\lambda^{(n)}\}$, then there is a neighborhood $\mathcal{O}$ of $(z_0,w_0)$ such that $\lambda_1=0$ for every $\lambda\in V\cap\mathcal{O}$. This implies $z_0=0$ and $\{0\}\times\mathbb{D}\subset V$. Hence $I_V\subset(z)$ and $I\subset(z^2)$. Then $w^n\in[I]_\alpha^\perp$, and it is easily verified
    $$[S_z^*,S_z]w^n=\frac{w^n}{\alpha+2},~n=0,1,\ldots,$$
    which again contracts to the essential normality of $S_z$.
\end{proof}
Combining Theorem \ref{thm:main1} and Theorem \ref{thm:distinguishedvariety}, we obtain the main theorem \ref{main-theorem}.

\section{K-homology for distinguished quotient modules}
~~~~Let $I$ be a distinguished ideal, then Theorem \ref{thm:main1} ensures the essential normality of $[I]_\alpha^\perp$ for sufficiently large $\alpha$. As usual, we use $\sigma_e([I]_\alpha^\perp)$ to denote the Taylor joint essential spectrum of the commuting pair $(S_z,S_w)$. Similarly to \cite[Proposition 2.5]{Ar2}, it can be verified that $C^*([I]_\alpha^\perp)=C^*\{Id_{[I]_\alpha^\perp}, S_z,S_w\}$ is irreducible, and therefore contains the ideal ${\cal K}([I]_\alpha^\perp)$ of all the compact operators. Similarly to \cite{WZ1,WZ2}, we have the following lemma.
\begin{prop}\label{prop:espectrum}
Let $I$ be an ideal of $\mathbb C[z,w]$, and $\alpha\geq-1$, such that $[I]_\alpha^\perp$ is essentially normal, then
$$\sigma_e([I]_\alpha^\perp)\subset Z(I)\cap\partial\mathbb{D}^2.$$
Moreover, if $I$ is distinguished, then $\sigma_e([I]_\alpha^\perp)=Z(I)\cap\partial\mathbb{D}^2.$
\end{prop}
\begin{proof}
By Spectral Mapping Theorem \cite{Cur,Tay2}, for $p\in I$ we have $$p(\sigma_e(S_z,S_w))=\sigma_e(p(S_z,S_w))=\sigma_e(S_p)=\{0\},$$
therefore $\sigma_e([I]_\alpha^\perp)\subset Z(I).$

To prove $\sigma_e([I]_\alpha^\perp)\subset\partial\mathbb{D}^2$, first we notice that
\begin{eqnarray}\label{eq:Berezin}
    &&Id_{[I]_\alpha^\perp}-S_zS_z^*-S_wS_w^*+S_{zw}S_{zw}^*\\
    &=&P_{[I]_\alpha^\perp}(\Id_{\mathcal{A}_\alpha(\mathbb{D}^2)}-M_zM_z^*-M_wM_w^*+M_{zw}M_{zw}^*)P_{[I]_\alpha^\perp}.\notag
\end{eqnarray}
It is not hard to see that $\mathcal{A}_\alpha(\mathbb{D}^2)$ has orthonormal basis $$\left\{\frac{1}{\Gamma(\alpha+2)}\sqrt{\frac{\Gamma(m+\alpha+2)}{\Gamma(m+1)}\frac{\Gamma(n+\alpha+2)}{\Gamma(n+1)}}z^mw^n:m,n\in\mathbb{N}\right\}.$$
In the case $\alpha>-1$, it is routine to verify
$$M_z^*(z^mw^n)=\frac{m}{m+\alpha+1}z^{m-1}w^n,$$
and
$$M_w^*(z^mw^n)=\frac{n}{n+\alpha+1}z^mw^{n-1}.$$
Then direct computation shows
\begin{eqnarray*}
    &&(\Id_{\mathcal{A}_\alpha(\mathbb{D}^2)}-M_zM_z^*-M_wM_w^*+M_{zw}M_{zw}^*)z^mw^n\\
    &=&\left(1-\frac{m}{m+\alpha+1}-\frac{n}{n+\alpha+1}+\frac{mn}{(m+\alpha+1)(n+\alpha+1)}\right)z^mw^n\\
    &=&\frac{(\alpha+1)^2}{(m+\alpha+1)(n+\alpha+1)}z^mw^n
\end{eqnarray*}
for $m,n>0$.

In the case $\alpha=-1$, it is obvious that
$$\Id_{\mathcal{A}_\alpha(\mathbb{D}^2)}-M_zM_z^*-M_wM_w^*+M_{zw}M_{zw}^*=1\otimes1.$$

In either cases, we have obtained the compactness of
$$\Id_{\mathcal{A}_\alpha(\mathbb{D}^2)}-M_zM_z^*-M_wM_w^*+M_{zw}M_{zw}^*.$$
Then since $C^*([I]_\alpha^\perp)$ is essentially commutative, it follows from (\ref{eq:Berezin}) that $\sigma_e([I]_\alpha^\perp)\subset\partial\mathbb{D}^2$.

In the case $I$ is distinguished, from \cite[Proposition 4.3]{Kne2} one observes
$$Z(I)\cap\mathbb{T}^2\subset\partial[Z(I)\cap\mathbb{D}^2].$$
Then similarly to the proof of \cite[Lemma 5.1]{WZ2} one can verify
$$Z(I)\cap\partial\mathbb{D}^2\subset\sigma_e([I]_\alpha^\perp),$$
hence completes the proof of the proposition.
\end{proof}

Let $I$ be a distinguished ideal of the polynomial ring $\mathbb C[z,w]$. By taking $\alpha\geq-1$ such that the quotient module $[I]_{\alpha}^\perp$ is essential normal,  we get an extension
$$
0\to\mathcal{K}\hookrightarrow C^*([I]_\alpha^\perp)\to C(Z(I)\cap \partial\mathbb D^2)\to 0,
$$
which yields a $K$-homology element in $K_1(Z(I)\cap\partial\mathbb D^2)$.

Now, let $p_1,\ldots,p_n$ be distinguished polynomials, such that $$Z(p_i)\cap Z(p_j)\cap\mathbb T^2=\emptyset,~i\not=j.$$
Set $p=\prod\limits_{i=1}^n p_i$. Taking $\alpha\geq-1$ such that $[p]_{\alpha}^\perp$ is essentially normal, by Proposition \ref{prop:enormal} it is not hard to see all the $[p_i]_{\alpha}^\perp$'s are essentially normal.
By Lemma \ref{lem:a-orthogonal} proved in the next section, $[p_i]_\alpha^\perp$ and $[p_j]_\alpha^\perp$ are asymptotically orthogonal for $i\not=j$, and hence by \cite[Theorem 3.3]{GWp1}, the subspace $[p_1]_\alpha^\perp+\cdots+[p_n]_\alpha^\perp$ is closed. Moreover, we have
$$[p]_\alpha^\perp=[p_1]_\alpha^\perp+\cdots+[p_n]_\alpha^\perp.$$
Now, let $e_{p}$ be the K-homology element defined by $[p]_{\alpha}^\perp$, and similarly $e_{p_i}$ the K-homology element defined by $[p_i]_{\alpha}^\perp$, $i=1,\cdots, n$. Since $[p_i]_{\alpha}^\perp$ is asymptotically orthogonal to $[p_j]_{\alpha}^\perp$, $i\not=j$, it is easy to see that
$$
e_p=e_{p_1}\oplus\cdots\oplus e_{p_n}.
$$
Next, we will show that if $p$ is a distinguished polynomial, and $\alpha\geq -1$ such that $[p]_{\alpha}^\perp$ is essentially normal, then the $K$-homology element $e_p$ is nontrivial. By \cite{GWk1}, it is sufficient to prove the following result.
\begin{thm}
Let $I$ be a distinguished ideal, and $\alpha\geq-1$ such that $[I]_\alpha^\perp$ is essentially normal, then the short exact sequence
$$
0\to \mathcal{K}\hookrightarrow C^*([I]_\alpha^\perp)\to C(Z(I)\cap \partial\mathbb D^2)\to 0
$$
is not splitting.
\end{thm}
\begin{proof}
By \cite[Lemma 5.5]{GWk1}, it suffices to find a Fredholm operator in $C^*([I]_\alpha^\perp)$ with nonzero index. Let $I=pL$ be the Beurling decomposition. Write $s=\dim\mathbb{C}[z,w]/L$, and take distinct points $z_1,\ldots,z_{s+1}\in\mathbb{D}$. For each $z_k$ choose $w_k$ satisfying $(z_k,w_k)\in Z(I)\cap\mathbb{D}^2$, and set
$$q(z,w)=\prod_{j=1}^{s+1}(z-z_k).$$
By Proposition \ref{prop:espectrum} and Spectral Mapping Theorem, $\sigma_e(S_q)\subset q(\mathbb{T}^2)$ which does not contain the point $0$, and therefore $S_q$ is Fredholm. Inequality (\ref{eq:injectivity}) gives
$$\dim\ker S_q\leq s<\infty.$$
On the other hand, each reproducing kernel $K_{(z_k,w_k)}\in\ker S_q^*$, and therefore
$$\dim\ker S_q^*\geq s+1>\dim\ker S_q,$$
which completes the proof.
\end{proof}

\section{Summation of distinguished quotient modules}
~~~~For a distinguished ideal $I\subset\mathbb{C}[z,w]$ which is not primary, we would like to reduce the weight index $\alpha$ in Theorem \ref{thm:main1}. Let $I=\bigcap_{i=1}^kI_i$ be the primary decomposition, then each $I_i$ is distinguished. Applying Theorem \ref{thm:main1} to each $I_i$, we obtain the essential normality of $[I_i]_\alpha^\perp$ for
$$\alpha\geq\max_{1\leq i\leq k}\max\{m_{I_i},n_{I_i}\}-2.$$
We shall prove the essential normality of $[I]_\alpha^\perp$ for these $\alpha$ by summing the $[I_i]_\alpha^\perp$'s up, under the assumption the $Z(I_i)$'s do not intersect each other at $\mathbb{T}^2$.

By \cite{GWp1}, if $H_1,H_2$ are subspaces of a Hilbert space $H$ such that $P_{H_1}P_{H_2}$ is compact, then $H_1$ and $H_2$ are said asymptotically orthogonal to each other. The following theorem from \cite[Theorem 3.3]{GWp1} provides a way of summation of essentially normal quotient modules.
\begin{thm}(Guo and Wang, 2007)\label{thm:a-orthogonal}
Let $\mathcal{N}_1$ and $\mathcal{N}_2$ be two essentially normal quotient modules, if they are asymptotically orthogonal, then $\mathcal{N}_1+\mathcal{N}_2$ is closed and essentially normal.
\end{thm}

To apply this theorem in our context, we need the following observation on the essential commutativity of $P_{[I]_\alpha^\perp}$ and multiplication operators.
\begin{lem}\label{lem:projectioncommutativity}
    Let $I\subset\mathbb{C}[z,w]$ be a distinguished ideal, $\alpha\geq-1$, such that $[I]_\alpha^\perp$ is essentially normal, then for every $p\in\mathbb{C}[z,w]$,
$[M_p^*,P_{[I]_\alpha^\perp}]$ is compact.
\end{lem}
\begin{proof}
    By Proposition \ref{prop:espectrum} and Spectral Mapping Theorem,
    $$\sigma_e(\Id_{[I]_\alpha^\perp}-S_wS_w^*)=\sigma_e(\Id_{[I]_\alpha^\perp}-S_w^*S_w)=\{0\},$$
    hence $S_w$ is essentially unitary.

    By assumption we have
    \begin{eqnarray*}
        0&\leq&P_{[I]_\alpha^\bot}(\Id_{\mathcal{A}_\alpha(\mathbb{D}^2)}-M_w^*M_w)P_{[I]_\alpha^\bot}\\
        &\leq&P_{[I]_\alpha^\bot}-P_{[I]_\alpha^\bot}M_w^*P_{[I]_\alpha^\bot}M_wP_{[I]_\alpha^\bot}\\
        &=&P_{[I]_\alpha^\bot}-S_w^*S_w\\
        &\in&\mathcal{K},
    \end{eqnarray*}
    which yields
    $$P_{[I]_\alpha^\bot}(\Id_{\mathcal{A}_\alpha(\mathbb{D}^2)}-M_w^*M_w)P_{[I]_\alpha^\bot}\in\mathcal{K}.$$
    Consequently from
    \begin{eqnarray*}
        &&P_{[I]_\alpha^\bot}M_w^*P_{[I]_\alpha}M_wP_{[I]_\alpha^\bot}\\
        &=&P_{[I]_\alpha^\bot}M_w^*(\Id_{\mathcal{A}_\alpha(\mathbb{D}^2)}-P_{[I]_\alpha^\perp})M_wP_{[I]_\alpha^\bot}\\
        &=&P_{[I]_\alpha^\bot}-S_w^*S_w-P_{[I]_\alpha^\bot}(\Id_{\mathcal{A}_\alpha(\mathbb{D}^2)}-M_w^*M_w)P_{[I]_\alpha^\bot}\in\mathcal{K},
    \end{eqnarray*}
    we get
    $$[M_w^*,P_{[I]_\alpha^\bot}]=P_{[I]_\alpha^\bot}M_w^*P_{[I]_\alpha}\in\mathcal{K}.$$
    Similarly $[M_z^*,P_{[I]_\alpha^\bot}]\in\mathcal{K}$, and the conclusion of the lemma follows from the fact that $\mathcal{K}$ is a closed subalgebra.
\end{proof}

\begin{lem}\label{lem:a-orthogonal}
Let $I_1$ and $I_2$ be ideals of $\mathbb C[z,w]$, where $I_1$ is distinguished, and $\alpha\geq1$. Set $\mathcal{N}_1=[I_1]_\alpha^\perp$ and $\mathcal{N}_2=[I_2]_\alpha^\perp$. If $\mathcal{N}_1$ is essentially normal and
$$
Z(I_1)\cap Z(I_2)\cap\partial\mathbb D^2=\emptyset,
$$
then $\mathcal{N}_1$ and $\mathcal{N}_2$ are asymptotically orthogonal.
\end{lem}
\begin{proof}
Denote by
$$\pi:B(\mathcal{N}_1)\to Q(\mathcal{N}_1)$$
the natural homomorphism to the Calkin algebra $Q(\mathcal{N}_1)=B(\mathcal{N}_1)/\mathcal{K}$. By Proposition \ref{prop:espectrum},
$$\sigma_e(S_z^{\mathcal{N}_1},S_w^{\mathcal{N}_1})=Z(I_1)\cap\mathbb T^2.$$

Let $p\in I_2$ be a polynomial such that $Z(p)=Z(I_2)$, then by assumption, $p(z,w)\neq0$ for all $(z,w)\in\sigma_e(S_z^{\mathcal{N}_1},S_w^{\mathcal{N}_1})$. Consequently $S_p^{\mathcal{N}_1}$ is Fredholm on $\mathcal{N}_1$. Therefore $P_{\mathcal{N}_1}P_{\mathcal{N}_2}$ is compact if and only if $(S_p^{\mathcal{N}_1})^*P_{\mathcal{N}_2}$ is compact, i.e. $\pi\left((S_p^{\mathcal{N}_1})^*P_{\mathcal{N}_2}\right)=0$. By Lemma \ref{lem:projectioncommutativity}, $[M_p^*,P_{\mathcal{N}_1}]$ is compact, and therefore
$$\pi\left((S_p^{\mathcal{N}_1})^*P_{\mathcal{N}_2}\right)=\pi(P_{\mathcal{N}_1}M_p^*P_{\mathcal{N}_1}P_{\mathcal{N}_2})=\pi(P_{\mathcal{N}_1}M_p^*P_{\mathcal{N}_2})=0,$$
which completes the proof of the lemma.
\end{proof}
\begin{cor}\label{cor:summation1}
    Suppose $I_1$ and $I_2\subset\mathbb{C}[z,w]$ are distinguished ideals, such that $Z(I_1+I_2)\cap\mathbb{T}^2=\emptyset$. Let $\alpha\geq-1$ such that both $[I_1]_\alpha^\perp$ and $[I_2]_\alpha^\perp$ are essentially normal, then $[I_1\cap I_2]_\alpha^\perp$ is also essentially normal.
\end{cor}
\begin{proof}
    Let $I_1=p_1L_1,~I_1=p_2L_2$ be their Beurling decompositions, then $Z(p_1)\subset Z(I_1)$ and $Z(p_2)\subset Z(I_2)$. Since $[I_1]_\alpha^\perp$ is essentially normal, and
    $$[I_1]_\alpha^\perp\ominus[p_1]_\alpha^\perp=[p_1]_\alpha\ominus[p_1L_1]_\alpha$$
    is finite dimensional, $[p_1]_\alpha^\perp$ is essentially normal. Similarly, $[p_2]_\alpha^\perp$ is essentially normal. By Lemma \ref{lem:a-orthogonal}, $[p_1]_\alpha^\perp$ and $[p_2]_\alpha^\perp$ are asymptotically orthogonal. Then by Theorem \ref{thm:a-orthogonal}, $[p_1]_\alpha^\perp+[p_2]_\alpha^\perp$ is essentially normal. Next we prove $[p_1p_2]_\alpha^\perp=[p_1]_\alpha^\perp+[p_2]_\alpha^\perp$.

    By assumption, $Z(p_1)\cap Z(p_2)\cap\partial\mathbb{D}^2=\emptyset$, hence the Taylor essential spectrum \cite{Cur} satisfies that $$\sigma_e(M_{p_1},M_{p_2})\subset\{(p_1(z,w),p_2(z,w)):z,w\in\mathbb{T}^2\}\subset\mathbb{C}^2\backslash\{(0,0)\},$$
    which gives the Fredholmness of $M_{p_1}M_{p_1}^*+M_{p_2}M_{p_2}^*$. It is clear that
    $$\ker(M_{p_1}M_{p_1}^*+M_{p_2}M_{p_2}^*)=\ker M_{p_1}^*\cap\ker M_{p_2}^*=[(p_1)+(p_2)]_\alpha^\perp.$$
    For $f\in[p_1]_\alpha\cap[p_2]_\alpha$, from
    $$[(p_1)+(p_2)]_\alpha=(\ker M_{p_1}M_{p_1}^*+M_{p_2}M_{p_2}^*)^\perp=\mathrm{ran}~(M_{p_1}M_{p_1}^*+M_{p_2}M_{p_2}^*),$$
    it follows that
    $$f\in\mathrm{ran}~(M_{p_1}M_{p_1}^*+M_{p_2}M_{p_2}^*).$$
    Hence there exists $g\in\mathcal{A}_\alpha(\mathbb{D}^2)$ such that
    $$(M_{p_1}M_{p_1}^*+M_{p_2}M_{p_2}^*)g=f\in[p_2]_\alpha.$$
    Consequently $M_{p_1}M_{p_1}^*g\in[p_2]_\alpha$. Set $h=P_{[p_2]_\alpha^\perp}M_{p_1}^*g$, then we have
    \begin{eqnarray*}
        S_{p_1}^{[p_2]_\alpha^\perp}h&=&P_{[p_2]_\alpha^\perp}M_{p_1}P_{[p_2]_\alpha^\perp}M_{p_1}^*g\\
        &=&P_{[p_2]_\alpha^\perp}M_{p_1}(P_{[p_2]_\alpha^\perp}+P_{[p_2]_\alpha})M_{p_1}^*g\\
        &=&P_{[p_2]_\alpha^\perp}M_{p_1}M_{p_1}^*g\\
        &=&0,
    \end{eqnarray*}
    which induces $h=0$ by Corollary \ref{lem:injectivity2}. Therefore $M_{p_1}^*g\in[p_2]_\alpha$ and consequently $M_{p_1}M_{p_1}^*g\in[p_1p_2]_\alpha$. Similarly $M_{p_2}M_{p_2}^*g\in[p_1p_2]_\alpha$, and it follows $f\in[p_1p_2]_\alpha$. Now we have proved that $$[p_1]_\alpha\cap[p_2]_\alpha\subset[p_1p_2]_\alpha,$$
    and since the inverse inclusion is obvious, we get
    $$[p_1]_\alpha\cap[p_2]_\alpha=[p_1p_2]_\alpha,$$
    which is equivalent to
    $$[p_1p_2]_\alpha^\perp=[p_1]_\alpha^\perp+[p_2]_\alpha^\perp.$$
    Therefore the the essential normality of $[p_1p_2]_\alpha^\perp$ follows from that of $[p_1]_\alpha^\perp+[p_2]_\alpha^\perp$. Finally, from
    \begin{eqnarray*}
        \dim_{\mathbb{C}}~p_1p_2\mathbb{C}[z,w]/(I_1\cap I_2)&\leq&\dim_{\mathbb{C}}~p_1p_2\mathbb{C}[z,w]/(p_1p_2L_1L_2)\\
        &=&\dim_{\mathbb{C}}~\mathbb{C}[z,w]/(L_1L_2)\\
        &<&\infty,
    \end{eqnarray*}
    we have
    $$\dim_{\mathbb{C}}~[I_1\cap I_2]_\alpha^\perp\ominus[p_1p_2]_\alpha^\perp=\dim_{\mathbb{C}}~[p_1p_2]_\alpha\ominus[I_1\cap I_2]_\alpha<\infty.$$
    Then the essential normality of $[I_1\cap I_2]_\alpha^\perp$ follows from the essential normality of $[p_1p_2]_\alpha^\perp$.
\end{proof}
Combining Theorem \ref{thm:main1} and Corollary \ref{cor:summation1}, the following result is immediate.
\begin{thm}\label{thm:main2}
    Let $I\subset\mathbb{C}[z,w]$ be a distinguished polynomial, $I=\bigcap_{i=1}^kI_i$ be the primary decomposition. Let $V_i=Z(I_i)\cap\mathbb{D}^2,~1\leq i\leq k$. Suppose that
    $$Z(I_i)\cap Z(I_j)\cap\mathbb{T}^2=\emptyset,~1\leq i<j\leq k.$$
    Then $[I]_\alpha^\perp$ is essentially normal for all integers
    $$\alpha\geq\max_{1\leq i\leq k}\max\{m_{V_i},n_{V_i}\}-2.$$
\end{thm}

\section{Non-distinguished examples}
~~~~In this section, we give some non-distinguished examples of essentially normal quotient modules.
\begin{exam}
    For $|a|<2$, consider the quotient Hardy module $[z-w+a]^\perp$. It is obvious to see $S_z=S_w-a~\Id$, hence
    $$[S_z^*,S_z]=[S_z^*,S_w-a~\Id]=[S_z^*,S_w],$$
    which is compact by \cite{Ya}. The compactness of $[S_w,S_w^*]$ follows in the same manner. Therefore $[z-w+a]^\perp$ is essentially normal.

    We make some remarks on this example. Write $q=z-w+a$. Firstly, if $a=0$, then $[z-w]^\perp$ is isomorphic to the Bergman module on the unit disc, and if $|a|\geq 2$, then $[z-w-a]^\perp=\{0\}$. Secondly, for $0<|a|<2$ we have
    $$
    Z(q)\cap \mathbb T^2\not=\emptyset,\quad Z(q)\cap(\mathbb T\times \mathbb D)\not=\emptyset,\quad\text{and } \quad Z(q)\cap(\mathbb D\times \mathbb T)\not=\emptyset.
    $$
    Finally, there is no polynomial in $[q]^\perp$, and therefore it is difficult to study its essential normality by direct calculation.
\end{exam}
It is known that for a quasi-homogeneous essentially normal quotient module \cite{WZ3}, each component of the variety intersects at most one of $\mathbb{T}^2,~\mathbb{T}\times\mathbb{D}$ or $\mathbb{D}\times\mathbb{T}$. The foregoing example shows that, for non-homogeneous quotient modules, this is not always the case.
\begin{cor}
For $|a|<2$, set $q=z+w-a$ and let $p$ be a distinguished polynomial such that $Z(q)\cap Z(p)\cap\mathbb T^2=\emptyset$, then the quotient Hardy module $[pq]^\perp$ is essentially normal.
\end{cor}
\begin{proof}
Since $Z(q)\cap Z(p)\cap \partial\mathbb D^2=\emptyset$, by Lemma \ref{lem:a-orthogonal}, $[q]^\perp$ is asymptotically orthogonal to $[I]^\perp$. By \cite[Lemma 2.6]{GWp1},
$$
[pq]^\perp=[q]^\perp+[q]^\perp,
$$
which is essentially normal.
\end{proof}

\textbf{Acknowledgement.} This work is supported by NSFC (No. 12271298, 12231005, 11871308, 12071253) and NSF of Shanghai (No. 21ZR1404200). The third author is partially supported by the Young Scholars Program of Shandong University. We would like to thank L. Chen (Shandong University) for many suggestions to make the paper more readable. We also thank D. Zheng and R. Yang for valuable discussions on this topic.

  \vskip3mm

  \noindent{Kunyu Guo, School of Mathematical Sciences, Fudan University, Shanghai, 200433, P. R. China, Email: kyguo@fudan.edu.cn}

  \noindent{Penghui Wang, School of Mathematics, Shandong University, Jinan 250100, Shandong, P. R. China, Email: phwang@sdu.edu.cn}

  \noindent{Chong Zhao, School of Mathematics, Shandong University, Jinan 250100, Shandong, P. R. China, Email: chong.zhao@sdu.edu.cn}
\end{document}